\documentstyle[11pt]{article}  
\onecolumn
\newtheorem{lem}{{\sc Lemma}}[section]              
\newtheorem{thm}{{\sc Theorem}}[section]            
\newtheorem{cor}{{\sc Corollary}}[section]          

\newcommand{\e}{\equiv}
\newcommand{\nqv}{\equiv \!\!\!\!\!\!\//~~}

\title{\Large\bf
  Triplets and Symmetries of Arithmetic mod $p^k$}
\author{ {\sc Nico F. Benschop} \\[1ex]
  {\it AmSpade Research  ~(n.benschop@chello.nl)}
  ~~The Netherlands}
\date{6 June 2000}

\begin{document}
\maketitle

\begin{abstract}
The finite ring $Z_k=Z$(+,~.) mod $p^k$ of residue arithmetic with
odd prime power modulus is analysed. The cyclic group of units $G_k$
in $Z_k$(.) has order $(p-1).p^{k-1}$, implying product structure
$G_k\e A_k.B_k$ with $|A_k|=p-1$ and $|B_k |=p^{k-1}$, the "core"
and "extension subgroup" of $G_k$ respectively. It is shown that
each subgroup $S \supset 1$ of core $A_k$ has zero sum, and $p$+1
generates subgroup $B_k$ of all $n\e 1$ mod $p$ in $G_k$.
The $p$-th power residues $n^p$ mod $p^k$ in $G_k$ form an order
$|G_k|/p$ subgroup $F_k$, with $|F_k|/|A_k|=p^{k-2}$, so $F_k$
properly contains core $A_k$ for $k \geq 3$.
By quadratic analysis (mod $p^3$) rather than linear analysis (mod $p^2$,
re: Hensel's lemma [5]~), the additive structure of subgroups $G_k$
and $F_k$ is derived. ~Successor function $n$+1 combines with the two
arithmetic $symmetries$ $-n$ and $n^{-1}$ to yield a $triplet$
structure in $G_k$ of three inverse pairs ($n_i,~n_i^{-1}$) with:
$n_i +1 \e -(n_{i+1})^{-1}$, ~indices mod 3, and $n_0.n_1.n_2 \e 1$
mod $p^k$. In case $n_0 \e n_1 \e n_2 \e n$ this reduces to the cubic
root solution $n+1 \e - n^{-1} \e - n^2$ mod $p^k ~(p$=1 mod 6). . .
The property of exponent $p$ distributing over a sum of core residues :
$(x+y)^p \e x+y \e x^p + y^p$ mod $p^k$ is employed to derive the known
$FLT$  inequality for integers. In other words, to a $FLT$ mod $p^k$
equivalence for $k$ digits correspond $p$-th power integers of $pk$
digits, and the $(p-1)k$ carries make the difference, representing
the sum of mixed-terms in the binomial expansion.
\end{abstract}

additive structure, triplet,
    cubic roots of unity,~ $carry$, ~Hensel, Fermat, $FST, FLT$.

{\bf MSC-class}: ~11D41

\section*{Introduction}

The commutative semigroup $Z_k$(.) of multiplication mod $p^k$
(prime $p>$2) has for all $k>$0 ~just two idempotents: $1^2 \e 1$
and $0^2 \e 0$, and is the disjoint union of the corresponding
maximal subsemigroups (~Archimedian components~[3], [4]~).
Namely the group $G_k$ of units ($n^i\e 1$ mod $p^k$ for some $i>$0)
which are all relative prime to $p$, and maximal ideal $N_k$ as
nilpotent subsemigroup of all $p^{k-1}$ multiples of $p$ ($n^i\e 0$
mod $p^k$ for some $i>$0). Order $|G_k|=(p-1)p^{k-1}$ has two
coprime factors, sothat $G_k\e A_kB_k$, with 'core' $|A_k|=p-1$
and 'extension group' $|B_k|=p^{k-1}$.
Residues of $n^p$ form a subgroup $F_k \subset G_k$ of order
$|F_k|=|G_k|/p$, to be analysed for its additive structure. Each
$n$ in core $A_k$ satisfies $n^p\e n$ mod $p^k$, a generalization
of Fermat's Small Theorem ($FST$) for $k>1$, denoted as $FST_k$.

{\bf Base $p$} number representation is used, which notation is
useful for computer experiments, as reported in tables 1,2. This
models residue arithmetic mod $p^k$ by considering only the $k$
less significant digits, and ignoring the more significant digits.
Congruence class [$n$] mod $p^k$ is represented by natural number
$n<p^k$, encoded in $k$ digits (base $p$). Class [$n$] consists of
all integers with the same least significant $k$ digits as $n$.

{\bf Define} the {\bf 0-extension} of residue $n$ mod $p^k$ as the
natural number $n < p^k$ with the same $k$-digit representation
(base $p$), and all more significant digits (at $p^m, ~m \geq k)$
set to 0.

Signed residue $-n$ is only a convenient notation for the complement
$p^k-n$ of $n$, which are both positive. $C[n]$ or $C_n$ is a cyclic
group of order $n$, such as $Z_k(+) \cong C[p^k]$. The units mod $p$
form a cyclic group $G_1=C_{p-1}$, and $G_k$ of order $(p-1).p^{k-1}$
is also cyclic for $k>$1 ~[1]. ~~Finite {\it semigroup structure} is
applied, and {\it digit analysis} of prime-base residue arithmetic,
to study the combination of (+) and (.) mod $p^k$, especially the
additive properties of multiplicative subgroups of ring $Z_k(+, .)~$.

Only elementary residue arithmetic, cyclic groups, and (associative)
function composition are used (thm3.2), to begin with the known
cyclic (one generator) nature of group $G_k$ of units mod $p^k$ [1].
Lemma 1.1 on the direct product structure of $G_k$,  and cor1.2 on
all $p$-th power residues mod $p^k$ as all extensions of those
mod $p^2$, are known in some form but are derived for completeness.
Lemma 1.4 on $B_k=(p+1)^*$, and further results are believed to be new.

The {\bf two symmetries} of residue arithmetic mod $p^k$, defined as
automorphisms of order 2, are {\bf complement} $-n$ under (+) with
factor $-1$, and {\bf inverse} $n^{-1}$ under (.) with exponent $-1$.\\
Their essential role in the triplet- structure (thm3.1) of this finite
ring is emphasized throughout. \\The main emphasis is on additive
analysis of multiplicative semigroup $Z$(.) mod $p^k$.\\
Concatenation will be used to indicate multiplication.

\begin{tabular}{ll}
{\bf Symbols} & and {\bf Definitions}~~~~~~(~odd prime $p$~)\\
\hline
$Z_k$(+,~.) & the finite ring of residue arithmetic mod $p^k$\\
$M_k$ & multiplication $Z_k$(.) mod $p^k$, semigroup ($k$-digit
      arithmetic base $p$)\\
$N_k$ & maximal ideal of $M_k:~n^i\e 0$ mod $p^k$~(some $i>$0),
     $|N_k|=p^{k-1}$\\
$n \in M_k$ & unique product ~$n = g^i.p^{k-j}$ mod $p^k$ ~($g^i
     \in G_j$ coprime to $p$)\\
0-extension X & of residue $x$ mod $p^k$:
     the smallest non-negative integer $X \e x$ mod $p^k$\\
(finite) extension $U$ & of $x$ mod $p^k$:
     any integer $U \e x$ mod $p^k$\\
$C_m$ or $C[m]$ & cyclic group of order $m$:
        ~~e.g. ~$Z_k(+) \cong C[p^k]$\\
$G_k\e A_k.B_k$ & group of units: all ~$n^i \e 1$ mod $p^k$
        ~(some $i>$0), ~~$|G_k|\e (p-1)p^{k-1}$\\
$A_k$ & ~~~~~{\bf core} of $G_k$, ~~$|A_k|=p-1~~~(n^p$=$n$
        mod $p^k$ for $n \in A_k$)\\
$B_k\e (p+1)^*$ & extension group of all $n\e 1$ mod $p$ ,
        ~~$|B_k|=p^{k-1}$\\
$F_k$ & subgroup of all $p$-th power residues in $G_k$ ,
        ~~$|F_k|=|G_k|/p$ \\
$A_k \subset F_k \subset G_k$ &
    proper inclusions only for $k \geq 3~~(A_2\e F_2 \subset G_2)$\\
$d(n)$ & core increment $A(n+1)-A(n)$ of core func'n
        $A(n)\e n^q,~q=|B_k|$\\
$FST_k$ & core $A_k$ extends $FST ~(n^p \e n$ mod $p$) to
        mod $p^{k>1}$ for $p-1$ residues\\
solution in core & $x^p+y^p \e z^p$ mod $p^k$ ~with~ $x,y,z$
        ~in core $A_k$.\\
period of $n \in G_k$ & order $|n^*|$ of subgroup generated by
        $n$ in $G_k(.)$\\
normation & divide $x^p+y^p \e z^p$ mod $p^k$ by one term (in $F_k$),
      yielding one term $\pm 1$\\
complement $-n$ & unique in $Z_k$(+) : ~$-n+n\e 0$ mod $p^k$\\
inverse $n^{-1}$ & unique in ~$G_k$(.) : ~$n^{-1}.~n\e 1$ mod $p^k$\\
1-complement $n"$ & unique in $Z_k$(+) : ~$n"+n\e -1$ mod $p^k$\\
inverse-pair & pair ($a, ~a^{-1}$) of inverses in $G_k$ \\
{\bf triplet} & 3 inv.pairs: ~$a+b^{-1}\e b+c^{-1}\e c+a^{-1}\e -1,
    ~(abc\e 1$ mod $p^k$)\\
triplet$^p$ & a triplet of three $p$-th power residues in  subgroup $F_k$ (thm3.1)\\
triplet$^p$ equiv'ce & one of the three equivalences of a triplet$^p$\\
symmetry mod $p^k$ & $-n$ and $n^{-1}$: ~order 2 automorphism of
        $Z_k(+)$ resp. $G_k(.)$\\
$EDS$ property & Exponent Distributes over a Sum:
        $(a+b)^p\e a^p+b^p$ mod $p^k$
\end{tabular}

\section{ Structure of the group $G_k$ of units }

\begin{lem}
~~~~$G_k ~\cong ~A'_k \times B'_k ~\cong ~C[p-1]~.~C[p^{k-1}]$

~~~~~~~~~~~~ and $M_k$ (mod $p^k$) has a sub-semigroup isomorphic to
    $M_1$ (mod $p$).
\end{lem}
\begin{proof}
~Cyclic group $G_k$ of {\it units} $n$ ($n^i\e 1$ for some $i>0$)
has order $(p-1)p^{k-1}$, namely $p^{k}$ minus $p^{k-1}$ multiples
of $p$. Then $G_k=A'_k \times B'_k$, the direct product of two
relative prime cycles, with corresponding subgroups $A_k$ and $B_k$,
sothat $G_k\e A_k.B_k$ where\\
{\bf extension group} $B_k=C[~p^{k-1}~]$ consists of all $p^{k-1}$
residues mod $p^k$ that are 1 mod $p$, \\
and ~{\bf core} $A_k=C[p-1]$, ~so $M_k$ contains sub-semigroup
$A_k \cup 0 \cong M_1$.
\end{proof}

{\bf Core $A_k$}, as $p-1$ cycle mod $p^k$, is Fermat's Small Theorem
$n^{p}\e n$ mod $p$ extended to $k>$1 for $p$ residues (including 0),
to be denoted as $FST_k$.\\ Recall that $n^{p-1} \e 1$ mod $p$ for
$n \nqv$0 mod $p$ ($FST$), then lem1.1 implies:

\begin{cor}~~With $|B|=p^{k-1}=q$ and $|A|=p-1$:\\ \hspace*{1cm}
Core $A_k=\{~n^q~\}$ mod $p^k~~(n=1.. p$-1) ~extends $FST$ for $k>$1,
 ~and: \\ \hspace*{1cm}
~~~~$B_k= \{n^{p-1}\}$ mod $p^k$ ~consists of all $p^{k-1}$
     residues 1 mod $p$ in $G_k$.
\end{cor}

Subgroup $F_k \e \{n^p\}$ mod $p^k$ of all $p$-th power residues
in $G_k$, with $F_k \supseteq A_k$ (only $F_2 \e A_2$) and order
$|F_k|=|G_k|/p=(p-1)p^{k-2}$, consists of {\bf all} $p^{k-2}$
extensions mod $p^k$ of the $p-1$ ~~$p$-th power residues in $G_2$,
which has order $(p-1)p$. Consequently we have:

\begin{cor}
  Each extension of $n^p$ mod $p^2$ ~(in $F_2$)
 is a $p$-th power residue (in $F_k$)
\end{cor}

{\bf Core generation}:
~The $p-1$ residues $n^q$ mod $p^k ~(q=p^{k-1})$
define core $A_k$ for 0$<n<p$.\\Cores $A_k$ for successive $k$
are produced as the $p$-th power of each $n_0<p$ recursively:\\
~$(n_0)^p \e n_1, ~(n_1)^p \e n_2,~(n_2)^p \e n_3$, etc.,
where $n_i$ has $i$+1 digits. In more detail:

\begin{lem}.\\ \hspace*{.5cm}
    The $p-1$ values $a_0<p$ define core $A_k$ by~
    $(a_0)^{p^{k-1}}=a_0+\sum_{i=1}^{k-1}a_ip^i$ ~(digits $a_i<p$).
\end{lem}
\begin{proof}
Let $a=a_0+mp<p^2$ be in core $A_2$, so $a^p\e a$ mod $p^2$.
Then $a^p= (a_0+mp)^p=a_0^p+p.a_0^{p-1}.mp\e a_0^p+mp^2$ mod $p^3$,
using $FST$. Clearly the second core digit, of weight $p$, is not
found this way as function of $a_0$, but requires actual computation
(unless $a\e p \pm 1$ as in lem1.3-4). It depends on the $carries$
produced in computing the $p$-th power of $a_0$. Recursively, each
next core digit can be found by computing the $p$-th power of a
core $A_k$ residue with $k$+1 digit precision; here core $A_k$
remains fixed since $a^p\e a$ mod $p^k$.
\end{proof}

Notice $(p^2 \pm 1)^p\e p^3 \pm 1$ mod $p^5$. Moreover, initial
$(p+1)^p\e p^2+1$ mod $p^3$ yields in general for $(p \pm 1)^{p^m}$
the next property:

\begin{lem} ~~$(p+1)^{p^m}\e p^{m+1}+1$ {\rm ~~mod} $p^{m+2}$\\
\hspace*{1cm}
      and: ~~~$(p-1)^{p^m}\e p^{m+1}-1$ {\rm ~~mod} $p^{m+2}$
\end{lem}

\begin{lem} 
~Extension group $B_k$ is generated by $p$+1 (mod $p^k$),
 with $|B_k|=p^{k-1}$, \\ \hspace*{1cm} and each subgroup
 $S \subseteq B_k$, ~$|S|=|B_k|/p^s$ has
  sum $\sum S \e |S|$ {\rm ~mod} $p^k$.
\end{lem}
\begin{proof}
 The period of $p+1$, which is the smallest $x$ with $(p+1)^x\e 1$
mod $p^k$, implies $m+1=k$ (re lemma 1.3). So $m=k-1$, yielding
period $p^{k-1}$. No smaller exponent generates 1 mod $p^k$ since
$|B_k|$ has only divisors $p^s$.

$B_k$ consists of all $p^{k-1}$ residues which are 1 mod $p$.
The order of each subgroup $S \subset B_k$ must divide $|B_k|$,
sothat $|S|=|B_k|/p^s$ ~($0 \leq s < k$) and~ $S=\{1+m.p^{s+1}\}
 ~(m=0~..~|S|-1)$.

Then ~$\sum S= |S|+p^{s+1}.|S|(|S|-1)/2$ ~mod $p^k$, ~where~
$p^{s+1}.|S|=p.|B_k|=p^k$,

sothat~ $\sum S = |S|= p^{k-1-s}$ mod $p^k$. ~~Hence no subgroup
of $B_k$ ~sums to 0 mod $p^k$.
\end{proof}

\begin{cor} 
~~For core $A_k \e g^*$: each unit $n \in G_k \e A_kB_k$ has the form:
\\ \hspace*{1cm}   $n \e g^i(p+1)^j$ mod $p^k$ for a unique pair of
     non-neg. exponents $i<|A_k|$ and $j<|B_k|$.
\end{cor}

Pair $(i,j)$ are the exponents in the core- and extension- component
of unit $n$.

\begin{thm} 
        Each subgroup $S \supset 1$ of core $A_k$ sums to ~0
        {\rm ~mod} $p^k~~(k>$0).
\end{thm}
\begin{proof}
    ~For {\it even} $|S|$: $-1$ in $S$ implies pairwise zero-sums.
In general: $c.S=S$ for all $c$ in $S$, and $c \sum S =\sum S$,
so~ $S.x=x$, writing $x$ for $\sum S$. Now for any $g$ in $G_k$:
$|S.g|=|S|$ sothat $|S.x|$=1 implies $x$ not in $G_k$, hence~
$x=g.p^e$ ~for some $g$ in $G_k$ and $0<e<k$ or $x=0 ~(e=k)$.
Then:

~~$S.x=S(g.p^e)=(S.g)p^e$ with $|S.g|=|S|$ if $e<k$.
~So~ $|S.x|$=1 ~yields~ $e$=$k$ ~and~ $x=\sum S$=0.
\end{proof}

Consider the {\bf normation} of an additive equivalence $a+b \e c$
mod $p^k$ in units group $G_k$, by multiplying all terms with the
inverse of one of these terms, for instance to yield rhs $-1$:

{\bf (1)} ~~~1-complement form: ~~$a+b \e -1$ mod $p^k$ in $G_k$
           ~~~(digitwise sum $p-1$, no carry).

For instance the well known $p$-th power residue equivalence:~
~$x^p+y^p \e z^p$ ~in $F_k$ yields: \\[1ex]
{\bf (2)} ~~~~normal form:~~~~~~~ $a^p+b^p \e -1$ mod $p^k$ ~in $G_k$,

with a special case (in core $A_k$) considered next.

\begin{figure}
\begin{picture}(100,230)(-180,0)
\setlength{\unitlength}{1.5mm}
\put( 0,50){Core ~A =  (43)* =  43 42 66 24 25 01 ~(mod $7^2$)}
\put( 5,45){ Cubic rootpair: 42 + 24 =  66 =  $-1$ }
\put(-30,35){\shortstack{
    Complement $C(n)=-n$ \\
     Inverse ~~~$I(n)=n^{-1}$\\
    Succesor $S(n)=n+1$ }}
\put( 4,21){\line(1,0){35}}
\put(21,12){\line(0,1){17}}
\put(13, 5){\vector(1, 2){17}}
\put(12,39){\vector(1,-2){18}}
\put(11,39){\vector(0,-1){35}}
\put(31, 4){\vector(0, 1){35}}
\put(12,39){\line(2,-1){30}}
\put(12, 4){\line(2, 1){30}}
\put( 9,30){I~~~~~~C~~~~~~~~C~~~~~~~~I}
\put(-30,10){\shortstack{
    {\large Symmetries:}\\
    $-n$ (diagonal) ~C\\
    $n^{-1}$ (vertical) ~I\\
    $-n^{-1}$ (horizontal) IC=CI }} \put(4,22){$-1$}
\put(40,20){\framebox(3.5,3.5){01}}  \put(18,40){+1}
\put(10,39){\framebox(3.5,3.5){42}}  \put(14,39){- - - - -$>$- - - -}
\put(30,39){\framebox(3.5,3.5){43}}  \put(35,40){$42+1=-(42)^{-1}$}
\put( 0,20){\framebox(3.5,3.5){66}}  \put(40,35){$42^3$=1 mod $7^2$}
\put(10, 0){\framebox(3.5,3.5){24}}  \put(14, 3){- - - - -$>$- - - -}
\put(30, 0){\framebox(3.5,3.5){25}}  \put(20, 0){+1}
\put( 7,40){$a$}      \put(33,37){$-a^{-1}$}
\put( 5, 1){$a^{-1}$} \put(35, 1){$-a$}
\end{picture}
\caption{{ ~~Core $A_2$ mod $7^2$ ($C_6$),
        ~~~Cubic roots $C_3$=\{42,~24,~01\}}}
  \label{fig.1}
\end{figure}
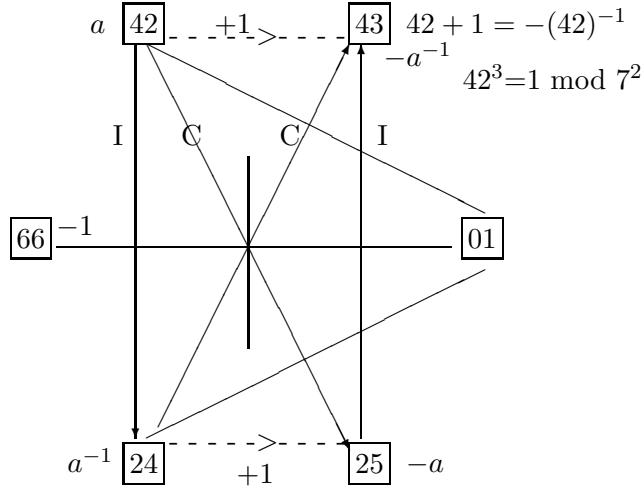

\section{The cubic root solution in core, ~and core symmetries}

\begin{lem} 
The cubic roots of 1 mod $p^k$ ($p \e 1$ mod 6) ~are $p$-th
  power residues in core $A_k$,
\\ \hspace*{1cm}
   and for $a^3 \e 1 ~(a \nqv 1):~ a+a^{-1} \e -1$ mod $p^{k>1}$
   has no 0-extension to integers.
\end{lem}
\begin{proof}
~If $p \e 1$ mod 6 then $3|(p-1)$ implies a core-subgroup
$S=\{a^2,a,1\}$ of three $p$-th powers: the cubic roots of 1
($a^3 \e 1$) in $G_k$ that sum to 0 mod $p^k$ (thm1.1). Now
$a^3-1=(a-1)(a^2+a+1)$, so if $a \nqv 1$ then $a^2+a+1 \e 0$,
hence $a+a^{-1} \e -1$ solves (1): a {\bf root-pair} of inverses,
with $a^2 \e a^{-1}$. ~$S$ in core consists of $p$-th power residues
with $n^{p}\e n$ mod $p^{k}$. Write $b$ for $a^{-1}$, then
$a^p+b^p \e -1$ {\bf and} $a+b \e -1$, sothat $a^p+b^p\e (a+b)^p$
mod $p^k$. ~Notice the {\bf \it Exponent Distributes over a Sum ($EDS$)},
implying inequality $A^p+B^p<(A+B)^p$ for the corresponding
0-extensions $A,~B,~A+B$ of core residues $a,~b,~a+b$ mod $p^k$.
\end{proof}

\begin{enumerate} {\small
\item
  {\bf ~Display} $G_k\e g^*$ by equidistant points on a {\bf unit circle}
  in the plane, with 1 and $-1$ on the horizontal axis (fig1,~2).
  The successive powers $g^i$ of generator $g$ produce $|G_k|$ points
  ($k$-digit residues) counter- clockwise. In this circle each inverse pair
  $(a,a^{-1})$ is connected $vertically$, complements $(a,-a)~diagonally$,
  and pairs $(a,-a^{-1})~horizontally$, representing functions $I, C$ and
  $IC=CI$ resp. (thm3.2). Figures 1, 2 depict for $p$=7, 5 these symmetries
  of residue arithmetic.
\item
  {\bf Scaling} any equation, such as $a+1\e -b^{-1}$, by a factor
  $s\e g^i \in G_k\e g^*$, yields $s(a+1)\e -s/b$ mod $p^k$,
represented by a {\bf rotation} counter clockwise over $i$ positions.
 }
\end{enumerate}

\subsection{Core increment symmetry at double precision, and asymmetry beyond}

Consider {\bf core function} $A_k(n)=n^{|B_k|} ~(~|B_k|=p^{k-1} ~cor1.1)$
as integer polynomial of odd degree, and {\bf core increment}
function $d_k(n)=A_k(n+1)-A_k(n)$ of even degree one less than $A_k(n)$.
Computing $A_k(n)$ upto precision $2k+1$ (base $p$) shows $d_k(n)$ mod
$p^{2k+1}$ to have a 'double precision' symmetry for 1-complements
$m+n=p-1$. Only $n<p$ need be considered due to periodicity $p$.\\
This naturally reflects in the additive properties of core $A_k$,
as in table 1 for $p$=7 and $k$=1, with $n<p$ in $A_1$ by $FST$:
symmetry of core increment $d_1(n)$ mod $p^3$ but not so mod $p^4$.
Due to $A_k(n) \e n$ mod $p ~(FST)$ we have $d_k(n) \e 1$ mod $p$,
so $d_k(n)$ is referred to as core 'increment', although in general
$d_k(n) \nqv 1$ mod $p^{k>1}$.

\begin{lem} ( Core increment at {\bf double precision} )~~
  For $q=|B_k|=p^{k-1}$ and $k>0$:\\
(a)~~ Core function $A_k(n) \e n^q$ mod $p^k$ and increment
   $d_k(n) \e A_k(n+1)-A_k(n)$ have {\bf period $p$}\\
(b)~~ for $m+n=p ~~~~:~ A_k(m) \e -A_k(n)$ mod $p^k$ {\rm ~~~(odd symmetry)}\\
(c)~~ for $m+n=p-1:~ d_k(m) \e d_k(n)$ mod $p^{2k+1}$
     and $\nqv$ mod $p^{2(k+1)}$ \\ \hspace*{1cm}
{\rm ('double precision' ~even symmetry and -inequivalence respectively).}
\end{lem}
\begin{proof}
{\bf(a)} ~Core function $A_K(n)\e n^q$ mod $p^k ~(q=p^{k-1},~n \nqv$ 0
mod $p$) has just $p-1$ distinct residues with $(n^q)^p \e n^q$ mod $p^k$,
and $A_k(n) \e n$ mod $p$ ($FST$). Including (non-core) $A_k(0) \e 0$ makes
$A_k(n)$ mod $p^k$ periodic in $n$ with {\bf period $p$} :
     ~$A_k(n+p) \e A_k(n)$ mod $p^k$, so $n<p$ suffices for core analysis.
Increment $d_k(n)$, as difference of two functions of period $p$,
also has period $p$.

{\bf(b)} ~$A_k(n)$ is a polynomial of odd degree with
  {\bf odd symmetry}~ $A_k(-n) \e (-n)^q \e -n^q \e -A_k(n)$.

{\bf(c)}~ Difference polynomial $d_k(n)$ is of even degree $q-1$ with
leading term $q.n^{q-1}$, and residues 1 mod $p$ in extension group
$B_k$. The even degree of $d_k(n)$ results in {\bf even symmetry},
because \\[1ex] \hspace*{1cm}
  $d_k(n-1) = n^q-(n-1)^q = -(-n)^q+(-n+1)^q = d_k(-n)$.

Denote $q=p^{k-1}$, then for ~$m+n=p-1$ follows:~
    $d_k(m)=A_k(m+1)-A_k(m)=(p-n)^q-m^q$
and $d_k(n)=A_k(n+1)-A_k(n)=(p-m)^q-n^q$,
yielding:~ $d_k(m)-d_k(n) = [~(p-n)^q+n^q~] -[~(p-m)^q+m^q~]$.
By binomial expansion and $n^{q-1} \e m^{q-1} \e 1$ mod $p^k$
in core $A_k$:~ $d_k(m)-d_k(n) \e 0$ mod $p^{2k+1}$. ~With~
$n \nqv m$ mod $p$:~ $m^{q-2} \e m^{-1} \nqv n^{-1} \e n^{q-2}$
mod $p$, ~causing~ $d_k(m) \nqv d_k(n)$ mod $p^{2k+2}$.
\end{proof}

Table 1 ~($p$=7) shows, for $\{n,m\}$ in core $A_1 ~(FST)$ with
$n+m \e -1$ mod $p$, the core increment symmetry mod $p^3$ and difference
mod $p^4 ~(k$=1). While $024^7 \e 024$ in $A_3 ~(k$=3) has core increment
1 mod $p^7$, but not 1 mod $p^8$, and similarly at 1-complementary cubic
root $642^7 \e 642$.

\subsection{ Another derivation of the cubic root of 1 mod $p^k$ }

The cubic root solution was derived, for 3 dividing $p-1$, via
subgroup $S \subset A_k$ of order 3 (thm1.1). For completeness a
derivation using elementary arithmetic follows.

Notice ~$a+b \e -1$ ~~to yield ~~~$a^2+b^2 \e (a+b)^2-2ab \e 1-2ab$,
~and:\\ \hspace*{1cm}
$a^3+b^3 \e (a+b)^3-3(a+b)ab \e -1+3ab$. ~~The combined sum is $ab-1$:
\\[1ex] \hspace*{.5cm}
$\sum_{i=1}^3(a^i+b^i) \e \sum_{i=1}^3 a^i + \sum_{i=1}^3 b^i \e ab-1$
~mod $p^k$. ~~Find $a,b$ for $ab \e 1$ mod$p^k$.

Since $n^2+n+1=(n^3-1)/(n-1)$=0 for $n^3 \e 1 ~(n \neq 1$), we have
$ab \e 1$ mod$p^{k>0}$ if $a^3 \e b^3 \e 1$ mod $p^k$, with 3 dividing
$p-1~(p \e 1$ mod 6). Cubic roots $a^3 \e 1$ mod $p^k$ exist for any
prime $p \e 1$ mod 6 at any precision $k>0$.

In the next section other solutions of $\sum_{i=1}^3 a^i +
\sum_{i=1}^3 b^i \e 0$ mod $p^k$ will be shown, depending not
only on $p$ but also on $k$, with $ab \e 1$ mod $p^2$
but $ab \nqv 1$ mod $p^3$, for some primes $p \geq 59$.

\section{Triplets, and the Core}

Any solution of (2): ~$a^p+b^p=-1$ mod $p^k$ has at least one term
($-1$) in core, and at most all three terms in core $A_k$.
To characterize such solution by the number of terms in core $A_k$,
quadratic analysis (mod $p^3$) is essential since proper inclusion
$A_k \subset F_k$ requires $k \geq 3$.
The cubic root solution, with one inverse pair (lem2.1), has all
three terms in core $A_{k>1}$. However, a computer search (table 2)
does reveal another type of solution of (2) mod $p^2$ for some
$p \geq 59$: three inverse pairs of $p$-th power residues,
denoted triplet$^p$, ~in core $A_2$.

\begin{thm} 
A {\bf triplet$^p$} of three inverse-pairs of $p$-th power residues
in $F_k$ satifies :
\\ \hspace*{1in} {\bf (3a)}~~~~~~$a+b^{-1} \e -1$ ~(mod $p^k$)
\\ \hspace*{1in} {\bf (3b)}~~~~~~$b+c^{-1} \e -1$ ~~~,,
\\ \hspace*{1in} {\bf (3c)}~~~~~~$c+a^{-1} \e -1$ ~~~,,
~~~with $abc \e 1$ mod $p^k$.
\end{thm}
\begin{proof} ~Multiplying by $b,~c,~a$ resp. maps
(3a) to (3b) if $ab \e c^{-1}$, and
(3b) to (3c) if $bc \e a^{-1}$, and
(3c) to (3a) if $ac \e b^{-1}$. ~All three conditions imply
     $abc \e 1$ mod $p^k$. \end{proof}

Table 2 shows all normed solutions of (2) mod $p^2$ for $p<200$, with
triplets at $p$= 59, 79, 83, 179, 193. The cubic roots, indicated by
$C_3$, occur only at $p \e 1$ mod 6, while a triplet$^p$ can occur for
either prime type $\pm 1$ mod 6. More than one triplet$^p$ can occur
per prime (two at $p$=59, three at 1093, four at 36847: each first
occurrance of such multiple triplet$^p$). There are primes for which
both rootforms occur, e.g. $p=79$ has a cubic root solution as well
as a triplet$^p$.

The question is if such {\bf loop structure} of inverse-pairs can have
a length beyond 3. Consider the successor $S(n)=n$+1 and the two
arithmetic symmetries, complement $C(n)=-n$ and inverse $I(n)=n^{-1}$,
as {\bf functions}, which compose associatively.\\ Then looplength $>$3
is impossible in arithmetic ring $Z_k(+,~.)$ mod $p^k$, seen as follows.

\begin{thm} .. (two basic solution types)\\  \hspace*{1cm}
    Each normed solution of (2) is (an extension of) a triplet$^p$
    or an inverse- pair.
\end{thm}
\begin{proof}
Assume $r$ equations $1-n_i^{-1}\e n_{i+1}$ form a loop of length $r$
(indices mod $r$). Consider function $ICS(n)\e 1-n^{-1}$, composed of
the three elementary functions: Inverse, Complement and Successor, in
that sequence.~
Let $E(n)\e n$ be the identity function, and $n \neq 0,1,-1$ to prevent
division by zero, then under function composition the {\bf third
iteration} $[ICS]^3=E$, since $[ICS]^2(n)\e -1/(n-1) ~\rightarrow
~[ICS]^3(n)\e n$ (repeat substituting $1-n^{-1}$ for $n$). Since $C$
and $I$ commute, $IC$=$CI$, the 3! =  6 permutations of \{$I,C,S$\}
yield only four distinct dual-folded-successor {\it "dfs"} functions:

\hspace*{.5cm} $ICS(n)=-n^{-1},~SCI(n)=-(1+n)^{-1},
              ~CSI(n)=(1-n)^{-1}, ~ISC(n)=-(1+n^{-1})$.

By inspection each of these has $[dfs]^3=E$, referred to as
{\bf loop length} 3. For a cubic rootpair {\it dfs=E}, and 2-loops
do not occur since there are no duplets (next note 3.2). Hence
solutions of (2) have only {\it dfs} function loops of length
1 and 3: inverse pair and triplet.
\end{proof}

A special triplet occurs if one of $a,b,c$ equals 1, say $a \e 1$.
Then $bc \e 1$ since $abc \e 1$, while (3a) and (3c) yield
$b^{-1} \e c \e -2$, so $b \e c^{-1} \e -2^{-1}$. Although triplet
$(a,b,c) \e (1,-2,-2^{-1})$ satisfies conditions (3), 2 is not in
core $A_{k>2}$, and by symmetry $a,b,c \nqv 1$ for any triplet$^p$
of form (3).~ If $2^p \nqv 2$ mod $p^2$ then 2 is not a $p$-th power
residue, so triplet $(1,-2,-2^{-1})$ is not a triplet$^p$ for such
primes (all upto at least $10^9$ except 1093, 3511).

\begin{figure} 
\begin{picture}(100,230)(-120,0)
\setlength{\unitlength}{1.5mm}
\put( 5,49){Core ~A ~~= ~~33 44 12 01 ~(mod $5^2$) }
\put(11,48){\framebox(3.5,3.5){}}
\put( 5,45){Extn ~B ~~=  11 21 31 41 01 ~~~~~~~~G =  A.B}
\put(13,46){\circle{4}}
\put(14,21){\line(1,0){14}}
\put(21,16){\line(0,1){10}}
\put(45,30){\vector(-1,2){3}}
\put(40,20){\framebox(3.5,3.5){01}}  \put(42,22){\circle{5}}
\put(39,27){03 =  g}
\put(37,33){14}
\put(33,37){02}
\put(27,39){11}  \put(28,40){\circle{5}}
\put(20,39){\framebox(3.5,3.5){33}}  \put(38,40){$33= -33^{-1}$}
\put(13,39){04}
\put( 7,37){22}
\put( 3,33){21}  \put( 4,34){\circle{5}}
\put( 1,27){13}  \put(14,27){+1~v~~~~~~~~~~~~~+1}
\put( 0,20){\framebox(3.5,3.5){44}}   \put(5,20){-1}
\put( 1,13){42}  \put(4,14){- - - -$>$ - - - - - triplet - - - - - -}
\put( 3, 7){31}  \put(4, 8){\circle{5}}
\put( 7, 3){43}
\put(13, 1){34}  \put(15,3){- -$>$- - - -}  \put(15,3){\line(1,6){6}}
\put(20,-1){\framebox(3.5,3.5){12}}
   \put(50,5){\shortstack{Triplet :\\(33,~34)\\(41,~42)\\(32,~33)}}
   \put(45,0){33.41.32 =  01}
\put(27, 1){41}  \put(28, 2){\circle{5}}
\put(33, 3){23}  \put(27, 5){\vector(-3,1){23}}
\put(37, 7){24}
\put(39,13){32}  \put(39,16){\vector(-3,4){17}}
\end{picture}
  \caption{{ ~~G =  A.B =  $g^*$ ~(mod $5^2$), ~~~~~Cycle in the plane}}
  \label{fig.2}
\end{figure}

\subsection{ A triplet for each $n$ in $G_k$ }

Notice the proof of thm3.2 does not require $p$-th power residues.
So {\bf any} $n \in G_k$ generates a triplet by iteration of one of
the four {\it dfs} functions (thm3.2), yielding the main triplet
structure of $G_k$ :

\begin{cor}
 ~{\rm Each} $n$ in $G_k ~(k>$0) generates a triplet of ~{\rm three}
  inverse pairs,\\
\hspace*{1cm} except if ~$n^3 \e 1$ and $n \nqv 1$ mod $p^k
       ~(p \e 1$ mod 6), which involves {\rm one} inverse pair.
\end{cor}
Starting at $n_0 \in G_k$ six triplet residues are generated  upon
iteration of e.g. $SCI(n)$: $n_{i+1}\e -(n_i+1)^{-1}$ (indices
mod 3), or another {\it dfs} function to prevent a non-invertable
residue. Less than 6 residues are involved if 3 or 4 divides $p-1$:

If $3|(p-1)$ then a cubic root of 1 ($a^3 \e 1, ~a \nqv 1$)
generates just 3 residues: ~$a+1\e -a^{-1}$; \\
--- together with its complement this yields a subgroup
   $(a+1)^*\e C_6$ ~(fig.1, $p$=7)\\
If 4 divides $p-1$ then an $x$ on the vertical axis has $x^2 \e -1$
so $x \e -x^{-1}$,\\
--- so the 3 inverse pairs involve then only five residues
~(fig.2: $p$=5).

\begin{enumerate} {\small
\item
   It is no coincidence that the period 3 of each {\it dfs} composition
[~of $-n,~n^{-1},~n$+1; \\ ~~~e.g: $CIS(n) \e 1-n^{-1}$~] exceeds the
number of symmetries of finite ring $Z_k(+,~.)$ by one.
\item
    {\bf No duplet} occurs: multiply $a+b^{-1} \e -1,~b+a^{-1} \e -1$
by $b$ resp. $a$ then $ab+1 \e -b$ and $ab+1 \e -a$, ~sothat ~$-b \e -a$
and $a \e b$.
\item
    {\bf Basic triplet} mod $3^2: G_2 \e 2^* \e \{2,4,8,7,5,1\}$ is a
6-cycle of residues mod 9. ~Iteration: ~$SCI(1)^*: -(1+1)^{-1} \e 4,~
-(4+1)^{-1} \e 7,~ -(7+1)^{-1} \e 1$, and $abc \e 1.4.7 \e 1$ mod 9. }
\end{enumerate}

\subsection{ The $EDS$ argument extended to non-core triplets }

The $EDS$ argument for the cubic root solution $CR$ (lem2.1), with all
three terms in core, also holds for any triplet$^p$ mod $p^2$. Because
$A_2 \e F_2$ mod $p^2$, so all three terms are in core for some linear
transform (5). Then for each of the three equivalences (3a-c) holds
the $EDS$ property: $(x+y)^p \e x^p+y^p$, and thus no finite (equality
preserving) extension exists, yielding inequality for the corresponding
integers for all $k>$1, to be shown next. A cubic root solution is a
special triplet$^p$ for $p \e 1$ mod 6, with $a \e b \e c$ in (3a-c).

Denote the $p-1$ core elements as residues of integer function
~$A(n)=n^{|B|}, ~(0<n<p)$, then by freedom of $p$-th power extension
beyond mod $p^2$ (cor1.4) choose, for any $k>2$ :

{\bf (4)} ~Core increment form:~ $A(n+1)-A(n) \e (r_n)^p$ mod $p^k$,

~~~~~~~~ with $(r_n)^p \e r_n.p^2+1 ~(r_n>$0),
~hence $(r_n)^p \e 1$ mod $p^2$, ~but ~$\nqv 1$ mod $p^3$ ~in general.

This rootform of triplets, with two terms in core, is useful for the
additive analysis of subgroup $F_k$ of $p$-th power residues mod $p^k$
(re: the known Fermat's Last Theorem $FLT$ case1: residues coprime
to $p$ - to be detailed in the next section).

Any assumed $FLT~case_1$ solution (5) can be transformed into form (4) in
two steps that preserve the assumed $FLT$ equality for integers $<p^{kp}$
in full $p$-th power precision $kp$ where $x,y<p^k$ , or $(k+1)p$
~in case $p^k<x+y<p^{k+1}$ (one carry).\\
Namely first $scaling$ by an integer $p$-th power factor $s^p$ that is
1 mod $p^2$ (so $s \e 1$ mod $p$), to yield as one lefthand term the core
residue $A(n+1)$ mod $p^k$. And secondly a $translation$ by an additive
integer term $t$ which is 0 mod $p^2$ applied to both sides, resulting
in the other lefthand term $-A(n)$ mod $p^k$, preserving the assumed
integer equality (unit $x^p$ has inverse $x^{-p}$ in $G_k$). Without
loss assume the normed form with $z^p \e 1$ mod $p^2$, then such
{\bf linear transformation} ($s,t$) yields:

{\bf (5)}~~~~~~~~~~~
 $x^p+y^p=z^p ~~\longleftrightarrow~~ (sx)^p+(sy)^p+t=(sz)^p+t$ ~~[ integers ],

~~~~~~~~ with~ $s^p \e A(n+1)/x^p, ~~~(sy)^p+t \e -A(n)$ ~mod $p^k$, ~so:

{\bf (5')} \hspace{3cm} $A(n+1)-A(n) \e (sz)^p+t$ ~mod $p^k$.

With $s^p \e z^p \e 1$ and $t \e 0$ mod $p^2$ this yields an equivalence
which is 1 mod $p^2$, hence a $p$-th power residue, with two of the three
terms in core. Such core increment form (4),(5') will be shown to have no
(equality preserving) finite extension, of all residues involved, to $p$-th
power integers, so the assumed integer $FLT~case_1$ equality cannot exist.

\begin{lem} 
 $p$-th powers of a 0-extended triplet$^p$ equivalence (mod $p^{k>1}$)
  yield integer inequality.
\end{lem}
\begin{proof}
In a triplet for some prime $p>2$ the core increment form (4) holds
for three distinct values of $n<p$, where scaling by respective factors
$-(r_n)^{-p}$  in $G_k$ mod $p^k$ returns 1-complement form (2).
Consider each triplet equivalence separately, and for a simple notation
let $r$ be any of the three $r_n$, with successive core residues
$A(n+1) \e x^p \e x, ~-A(n) \e y^p \e y$ mod $p^k$.
Then~ $x^p+y^p \e x+y \e r^p$ mod $p^k$, where $r^p \e 1$ mod $p^2$, has
both summands in core, but right hand side $r^p \nqv 1$ mod $p^{k>2}$ is
not in core, with deviation $d \e r-r^p \nqv 0$  mod $p^k$.\\
Hence~ $r \e r^p+d \e (x+y)+d$ mod $p^k$ ~(with $d \e 0$ mod $p^k$
in the cubic root case), and~ $x^p+y^p \e (x+y+d)^p$ mod $p^k$.
This equivalence has no finite (equality preserving) 0-extension to
integer $p$-th powers since $X^p+Y^p<(X+Y+D)^p$, ~so the assumed $FLT$
case$_1$ solution cannot exist.
\end{proof}
For $p$=7 the cubic roots are $\{42,24,01\}$ mod $7^2$ (base 7).
In full 14 digits: $42^7+24^7=01424062500666$ while $66^7=60262046400666$,
which are equivalent mod $7^5$ but differ mod $7^6$.

More specifically, linear transform (5) adjoins to a $FLT$ case$_1$
solution mod $p^2$ a solution with two adjacent core residues
mod $p^k$ ~(5') for any precision $k>1$, while preserving the
assumed integer $FLT$ case$_1$ equality. Without loss one can
assume  scalefactor $s<p^k$ and  shift term $t<p^{2k}$, yielding
double precision integer operands $\{sx,sy,sz\} < p^{2k}$, with
an (assumed) $p$-th power equality of terms $<p^{2kp}$.
  Although equivalence mod $p^{2k+1}$ can hold by proper choice
  of linear transform $(s,t)$,  ~$inequivalence$ at base $p$
  ~triple precision~ $3k+1$ ~follows by:

\begin{lem}  ~(~triple precision inequality ~): \\ \hspace*{5mm}
Any extension $(X,Y,Z)$ of $(x,y,z)$ in~ $x^p+y^p \e z^p$ mod $p^{k>1}$
($FLT_k$ case$_1$) yields an integer $p$-th power inequality
(of terms $<p^{pk}$), with in fact inequivalence $X^p+Y^p \nqv Z^p$
   mod $p^{3k+1}$.
\end{lem}
\begin{proof} 
Let $X=up^k+x,~ Y=vp^k+y,~ Z=wp^k+z$ extend the residues $x,y,z<p^k$,
~such that $X,Y$ mod $p^{k+1}$ are not both in core $A_{k+1}$.
So the extensions do not extend core precision $k$, and without
loss take $u,v,w<p^k$, due to a scalefactor $s<p^k$ in (5). Write
$h=(p-1)/2$, ~then binomial expansion upto quadratic terms yields:

~~~~~~~~ $X^p \e u^2x^{p-2}h~p^{2k+1}+ux^{p-1}p^{k+1}+x^p$ ~~mod $p^{3k+1}$,
~and similarly:

~~~~~~~~ $Y^p \e v^2y^{p-2}h~p^{2k+1}+vy^{p-1}p^{k+1}+y^p$ ~~mod $p^{3k+1}$,
   ~and:

~~~~~~~~ $Z^p \e w^2z^{p-2}h~p^{2k+1}+wz^{p-1}p^{k+1}+z^p$ ~~mod $p^{3k+1}$,

where: $x^p+y^p \e x+y \e z^p$ ~mod $p^k$,
~and~ $x^{p-1} \e y^{p-1} \e 1$ ~mod $p^k$, but not so mod $p^{k+1}$:
\\ ~$u,v$ ~are such that not both $X,Y$ are in core $A_{k+1}$, hence
core precision $k$ is not increased.

By lemma 3.1 the 0-extension of $x,y,z$ (so $u=v=w=0$) does not
yield the required equality $X^p+Y^p=Z^p$. To find for which
maximum precision equivalence $can$ hold, choose $u,v,w$ sothat:

~~$(u+v)p^{k+1}+x^p+y^p \e wz^{p-1}p^{k+1}+z^p$ mod $p^{2k+1}$ ...[*]..
~yielding~ $X^p+Y^p \e Z^p$ mod $p^{2k+1}$.

A cubic root solution has also $z^p \e z$ in core $A_k$, so
$z^{p-1} \e 1$ mod $p^k$, then $w=u+v$ with $w^2>u^2+v^2$ would
require $x^p+y^p \e z^p$ mod $p^{2k+1}$, readily verified for $k$=2
and any prime $p>2$. \\[1ex] Such extension [*] implies inequivalence
$X^p+Y^p \nqv Z^p$ mod $p^{3k+1}$ for non-zero extensions $u,v,w$.
Because $u+v=w$ together with $u^2+v^2=w^2=(u+v)^2$ yields $uv=0$.
So any (zero- or nonzero-) extension yields inequivalence mod $p^{3k+1}$.
\end{proof}

\section{Residue triplets and Fermat's integer powersum inequality}

Core $A_k$ as $FST$ extension, the additive zero-sum property
of its subgroups (thm1.1), and the triplet structure of units
group $G_k$, allow a direct approach to Fermat's Last Theorem:

{\bf (6)}~~~~~
        $x^p+y^p = z^p$ (prime $p>2$) ~has no solution for positive
     integers $x,~y,~z$\\ \hspace*{2cm}
     with ~~case$_1$ : ~$xyz \nqv 0$ mod $p$,
     and ~~case$_2$ : ~$p$ divides one of $x,y,z$.

Usually (6) mentions exponent $n>2$, but it suffices to show
inequality for primes $p>2$, because for composite exponent $m=p.q$
holds $a^{pq}=(a^p)^q= (a^q)^p$. If $p$ divides two terms then it
also divides the third, and all terms can be divided by $p^p$.
So in $case$ 2:~ $p$ divides just one term.

A finite integer $FLT$ solution of (6) has three $p$-th powers $<p^k$
for some finite fixed $k$, so occurs in $Z_k$, yet with no $carry$
beyond $p^{k-1}$, and (6) is the 0-extension of this solution mod
$p^k$. Each residue $n$ mod $p^k$ is represented uniquely by $k$
digits, and is the product of a $j$-digit number as 'mantissa'
relative prime to $p$, and $p^{k-j}$ represented by $k-j$ trailing
zero's (cor1.3).

Normation (2) to $rhs=-1$ simplifies the analysis, and maps residues
($k$ digits) to residues, keeping the problem finite. Inverse normation
back to (6) mod $p^k$ is in case 1 always possible, using an inverse
scale factor in group $F_k$. So normation does {\bf not} map to the
reals or rationals.

The present approach needs only a simple form of Hensel's lemma [5]
(in the general $p$-adic number theory), which is a direct consequence
of cor1.2 : ~~extend digit-wise the 1-complement form  such that the
$i$-th digit of weight $p^i$ in $a^p$ and $b^p$ sum to $p-1$
~(all $i \geq 0$), with $p$ choices per extra digit. Thus to each
normed solution of (2) mod $p^2$ correspond $p^{k-2}$ solutions
mod $p^k$:

\begin{cor} 
~(1-cmpl extension)~~~~
        A normed $FLT_k$ root is an extended $FLT_2$ root.
\end{cor}

\subsection{Proof of the FLT inequality}  

Regarding $FLT$ case$_1$, an inverse-pair and triplet$^p$ are the only
(normed) $FLT_k$ roots (thm3.2). As shown (lem3.1), any assumed integer
case$_1$ solution has a corresponding equivalent core increment form
(4) with two terms in core, having no integer extension, against the
assumption.

\begin{thm} ($FLT$ Case 1).
   ~For prime $p>$2 and integers $x,y,z>0$ coprime to $p$ :\\ \hspace*{1in}
   ~$x^p+y^p=z^p$ has no solution.
\end{thm}
\begin{proof}
~An $FLT_{k>1}$ solution is a linear transformed extension of an $FLT_2$
root in core $A_2=F_2$ (cor4.1). By lemmas 2.2c and 3.1 it has no finite
$p$-th power extension, yielding the theorem.
\end{proof}

In $FLT$ case$_2$ just one of $x,y,z$ ~is a multiple of $p$, hence
$p^p$ divides one of the three $p$-th powers in $x^p+y^p=z^p$. Again,
any assumed case$_2$ equality can be scaled and translated to yield
an equivalence mod $p^p$ with two terms in core $A_p$, having no
integer extension, contra the assumption.

\begin{thm} ($FLT$ case$_2$)
~~~For prime $p>$2 and positive integers $x,y,z$ :\\ \hspace*{1in}
 if~ $p$ ~divides just one of~ $x,y,z$ ~then ~$x^p+y^p=z^p$
 ~has no solution.
\end{thm}
\begin{proof}
In a case$_2$ solution $p$ divides a lefthand term, $x=cp$ or
$y=cp~(c>$0), or the right hand side $z=cp$. Bring the multiple of
$p$ to the right hand side, for instance if $y=cp$ we have~
$z^p-x^p=(cp)^p$, while otherwise $x^p+y^p=(cp)^p$. So the sum or
difference of two $p$-th powers coprime to $p$ must be shown not to
yield a $p$-th power $(cp)^p$ for any $c>0$:

{\bf (7)}~~~~
    $x^p \pm y^p = (cp)^p$ ~has no solution for integers $x,y,c>0$.

Notice that core increment form (4) does not apply here. However,
by $FST$ the two lefthand terms, coprime to $p$, are either
complementary or equivalent mod $p$, depending on their sum or
difference being $(cp)^p$. ~Scaling by $s^p$ for some $s \e 1$
mod $p$, ~so $s^p \e 1$ mod $p^2$, transforms one lefthand term into
a core residue $A(n)$ mod $p^p$, with $n \e x$ mod $p$.
 And translation by adding $t \e 0$ mod $p^2$ yields the other term
$A(n)$ or $-A(n)$ mod $p^p$ respectively. The right hand side then
becomes $s^p(cp)^p+t$, ~equivalent to $t$ mod $p^p$. So an
assumed equality (7) yields, by two equality preserving tansformations,
the next equivalence (8), where $A(n) \e u \e u^p$ mod $p^p$ ~($u$ in
core $A=A_p$ for $0<n<p$ with $n \e x \e u$ mod $p$) and $s \e 1,~ t \e 0$
mod $p^2$:

{\bf (8)}~~~~~ $u^p \pm u^p \e u \pm u \e t$ mod $p^p ~(u \in A_p)$,~
   ~where~ $u \e (sx)^p$ ,~ $\pm u \e \pm (sy)^p+t$ mod $p^p$.

Equivalence (8) does not extend to integers, because $U^p+U^p>U+U$,
~and $U^p-U^p=0 \neq T$, ~where $U,T$ are the 0-extensions of $u,t$
mod $p^p$ respectively. But this contradicts assumed equalities (7),
which consequently must be false. \end{proof}

{\bf Remark}: {\small
~~From a practical point of view the $FLT$ integer inequality of a
0-extended $FLT_k$ root (case$_1$) is caused by the {\it carries}
beyond $p^{k-1}$, amounting to a multiple of the modulus, produced in
the arithmetic (base $p$).
In the expansion of $(a+b)^p$, the mixed terms $can$ vanish mod $p^k$
for some $a,b,p$. Ignoring the carries yields $(a+b)^p \e a^p+b^p$
mod $p^k$,  and the $EDS$' property is as it were the $syntactical$
expression of ignoring the carry ($overflow$) in residue arithmetic.
In other words, in terms of $p$-adic number theory, this means
'breaking the Hensel lift': the residue  equivalence of an $FLT_k$
root mod $p^k$, although it holds for all $k>$0, $does$ imply
inequality for integers due to its special triplet structure,
where exponent $p$ distributes over a sum.}

\section*{ Conclusions }

\begin{enumerate}
\item
   Symmetries $-n,~n^{-1}$ determine $FLT_k$ roots but do not exist
    for positive integers.
\item
  Another proof of $FLT$ case$_1$ might use product 1 mod $p^{k}$ of
$FLT_k$ root terms: $ab \e 1$ or $abc \e 1$, which is impossible for
integers $>1$. The product of $m$ (=2, 3, $p$) ~$k$-digit integers
has $mk$ digits. ~Arithmetic mod $p^k$ {\bf ignores carries} of
weight $p^k$ and beyond.
Removal of the mod $p^k$ condition from a particular $FLT_k$ root
equivalence 0-extends its terms, and the ignored carries imply
inequality for integers.
\item
{\bf Core} $A_k \subset G_k$ as extension of $FST$ to mod $p^{k>1}$,
and the zero-sum of its subgroups (thm1.1) yielding the cubic $FLT$
root (lem2.1), started this work. The triplets were found by analysing
a computer listing (tab.2) of the $FLT$ roots mod $p^2$ for $p<200$.
\item
Linear analysis (mod $p^2$) suffices for root existence (Hensel,
cor4.1), but {\bf quadratic} analysis (mod $p^3$) is necessary to
derive triplet$^p$ core-increment form {\bf (4) ~(5,5')} with
maximally two terms in core $A_3$.
\item
"$FLT$ eqn(1) has no finite solution" and "$[ICS]^3$ has no finite
fixed point" \\are equivalent (thm3.2), yet each $n \in G_k$ is a
fixed point of $[ICS]^3$ mod $p^k$ \\ (re: $FLT_2$ roots imply all
roots for $k>$2, yet no 0-extension to integers).
\item
Crucial in finding the arithmetic triplet structure, and the double
precision core-increment symmetry and inequivalence (lem2.2c) were
extensive computer experiments, and the application of {\it associative
function composition}, the essence of semi-groups, to the three
elementary functions (thm3.2): \\ \hspace*{1cm}
  successor $S(n)=n$+1, complement $C(n)=-n$ and inverse $I(n)=n^{-1}$,
\\ with period 3 for $SCI(n)=-(n+1)^{-1}$ and the other three such
compositions. In this sense $FLT$ is not a purely arithmetic problem,
but essentially requires non-commutative and associative function
composition for its proof.
\end{enumerate}

\section*{ Acknowledgements }

The opportunity given me by the program committee in Prague [2],
to present this simple application of finite semigroup structure to
arithmetic, is remarkable and greatly appreciated. Also, the feedback
from several correspondents is gratefully acknowledged.

\section*{ References }

\begin{enumerate} {\small
\item T.Apostol: {\it Introduction to Analytical Number Theory}
(thm 10.4-6),  Springer Verlag, 1976.
\item N.F.Benschop: "The semigroup of multiplication mod $p^k$, an
   extension of Fermat's Small Theorem, and its additive structure",\\
  International conference {\it Semigroups and their Applications}
  (Digest p7)   Prague, July 1996.
\item A.Clifford, G.Preston: {\it The Algebraic Theory of Semigroups}
   \\ Vol 1 (p130-135), AMS survey \#7, 1961.
\item S.Schwarz: "The Role of Semigroups in the Elementary Theory of
   Numbers", \\ Math.Slovaca V31, N4, p369-395, 1981.
\item G.Hardy, E.Wright: {\it An Introduction to the Theory of Numbers}
 \\ (Chap 8.3, Thm 123), Oxford-Univ. Press 1979.
 }
\end{enumerate}
\begin{center} -----///----- \end{center}
\newpage
\begin{verbatim}

  n.   n     F= n^7      F'= PDo        PD1         PD2       p=7
  0.  0000   000000000   000000001   010000000   000000000  7-ary code
  1.  0001   000000001   000000241   023553100   050301000  9 digits
  2.  0002   000000242   000006001 < 055440100   446621000
  3.  0003   000006243   000056251   150660100   401161000   '<' :
  4.  0004   000065524   000345001 < 324333100   302541000  Cubic roots
  5.  0005   000443525   001500241   612621100   545561000 (n+1)^p - n^p
  6.  0006   002244066   004422601   355655100   233411000  = 1 mod p^3
         x          xx         ^^^sym
  7.  0010   010000000   013553101   410000000   000000000
  8.  0011   023553101   031554241   116312100   062461000
  9.  0012   055440342   062226001 < 351003100   534051000
 10.  0013   150666343   143432251   630552100   600521000
 11.  0014   324431624   255633001 < 455101100   160521000
 12.  0015   613364625   444534241   156135100   242641000
 13.  0016   361232166   025434501   110316100   223621000

 14.  0020   420000000   423165201   010000000   000000000
 15.  0021   143165201   263245241   402261100   313151000
 16.  0022   436443442   342105001 < 502606100   060611000
 17.  0023   111551443   000651251   326354100   541031000
 18.  .024   112533024 ! 660000001 < 036146100   035011000 (n+1)^p - n^p
 19.  .025   102533025 ! 366015241   612322100   531201000  = 1 mod p^7
 20.  0026   501551266   625115401   332500100   600441000
--------&c
Table 1:  Periodic Difference of i-th digit:  PDi(n) = F(n+p^i) - F(n)
\end{verbatim}
\newpage \label{lastpage}
\begin{verbatim}

Find a+b = -1 mod p^2 (in A=F < G): Core A={n^p=n}, F={n^p} =A if k=2.
 G(p^2)=g*, log-code: log(a)=i, log(b)=j;  a.b=1 --> i+j=0 (mod p-1)

TRIPLET^p: a+ 1/b= b+ 1/c= c+ 1/a=-1; a.b.c=1; (p= 59 79 83 179 193 ...
^^^^^^^
Root-Pair:  a+ 1/a=-1; a^3=1 ('C3') <--> p=6m+1 (Cubic rootpair of 1)
^^^^^^^^^
p:6m+-1 g=generator;   p < 2000:  two  triplets at p= 59, 701, 1811
  5:-   2                        three triplets at p= 1093
  7:+   3  C3    11:-  2
 13:+   2  C3    17:-  3
 19:+   2  C3    23:-  5   29:-   2
 31:+   3  C3
 37:+   2  C3    41:-  6
 43:+   3  C3    47:-  5
 53:-   2                   log    lin mod p^2
 59:-   2                 ------  ------------
   -2,-25( 40 15, 18 43)  25, 23( 35 11, 23 47) -23,  2( 53 54,  5  4)
              --     --              --     --              --     --
   27, 19( 18 44, 40 14) -19,  8( 13 38, 45 20)  -8,-27(  5  3, 53 55)
 61:+   2  C3
 67:+   2  C3    71:-  7
 73:+   5  C3
 79:+   3  C3
   30, 20( 40 46, 38 32) -20, 10( 36 42, 42 36) -10,-30( 77 11,  1 67)
 83:-   2
   21,  3(  9 74, 73  8)  -3, 18( 54 52, 28 30) -18,-21( 13 36, 69 46)
 89:-   3
 97:+   5  C3   101:-  2
103:+   5  C3   107:-  2
109:+   6  C3   113:-  3
127:+   3  C3   131:-  2   137:-  3
139:+   2  C3   149:-  2
151:+   6  C3
157:+   5  C3
163:+   2  C3   167:-  5   173:-  2
179:-   2
 19,  1( 78 176,100  2)  -1, 18( 64 90,114 88) -18,-19( 88 59, 90 119)
181:+   2  C3   191:- 19
193:+   5  C3
 -81, 58( 64 106,128 86) -58, 53( 4 101,188 91) -53, 81(188 70, 4 122)
197:-   2
199:+   3  C3
-------       -------------------------------------
Table 2:      FLT_2 root: inv-pair (C3) & triplet^p   (for p < 200)
\end{verbatim}
\end{document}